%% file: kolbds-short.tex
\documentclass[10pt]{article}
\include{macros}

\usepackage[all]{xy}
\newcommand{\np}[1]{#1}   
\renewcommand{\np}[1]{}

\usepackage[hypertex]{hyperref}

\title{Global Divisibility of Heegner Points and Tamagawa Numbers} 

\author{Dimitar P. Jetchev\footnote{The author was supported by graduate fellowships from Microsoft Research and the University of California at Berkeley.}
\\ 
\\ {\small University of California at Berkeley}
\\ {\small \tt jetchev@math.berkeley.edu}
}
\date{}

\begin{document}
\maketitle

\begin{abstract}   
We improve Kolyvagin's upper bound on the order of the $p$-primary part of the Shafarevich-Tate group 
of an elliptic curve of rank one over a quadratic imaginary field. In many cases, our bound is precisely the one predicted 
by the Birch and Swinnerton-Dyer conjectural formula. 
\end{abstract}


\section{Introduction}\label{sec:intro}

Let $E_{/\Q}$ be an elliptic curve of conductor $N$ and let $D > 0$ be a fundamental discriminant, such that all prime factors of $N$ 
are split in the quadratic imaginary field $K = \Q(\sqrt{-D})$. We call such a $D$ a \emph{Heegner discriminant} for the elliptic curve 
$E_{/\Q}$. 
   
Let $\cN \subset \cO_K$ be an ideal in the ring of integers of $K$, such that $N = \cN \cdot \overline{\cN}$ 
and $\cO_K / \cN \isom \Z / N\Z$. 
The Heegner point $x_1 = [\C/\cO_K \ra \C / \cN^{-1}]$ lies on the modular curve $X_0(N)$ and is defined over the Hilbert 
class field $H/K$ by the theory of complex multiplication. If $\vphi : X_0(N) \ra E$ is a fixed optimal modular 
parametrization which maps the cusp $i\infty$ of $X_0(N)$ to the origin of $E$ (see~\cite{breuil-conrad-diamond-taylor} for the existence of 
such a parametrization), 
then the Gross-Zagier formula (see~\cite[Ch.I,6.5]{gross-zagier} relates the height of the point 
$y_K = \Tr_{H/K}(\vphi(x_1))$ to the special value $L'(E_{/K}, 1)$ of the derivative of the $L$-function 
$L(E_{/K}, s)$. Gross and Zagier used this formula (see~\cite[Ch.V,2.2]{gross-zagier}) to restate the Birch and 
Swinnerton-Dyer conjectural formula for $E_{/K}$ (whenever the analytic rank is one) as follows  

\begin{conjecture}[Birch and Swinnerton-Dyer formula]\label{conj:GZ-BSD}
If the point $y_K$ has infinite order then $E(K)$ has rank one, the Shafarevich-Tate group 
$\Sha(E/K)$ is finite and its order is given by 
$$
\# \Sha(E/K) = \left ( \frac{[E(K) : \Z y_K]}{c \cdot \prod_{q \mid N}c_q} \right )^2,  
$$
where $c$ is the Manin constant (known to be a positive integer) and
$c_q = [E(\Q_q) : E^0(\Q_q)]$ is the Tamagawa number at $q$. 
\end{conjecture}

\noindent Kolyvagin (see~\cite[Thm.A]{kolyvagin:euler_systems},~\cite{kolyvagin:structure_of_selmer} 
and~\cite{kolyvagin:structureofsha}) has shown the rank part of the above conjecture, the finiteness 
of $\Sha(E/K)$ and a significant part of the conjectural formula. More precisely, 
consider the following hypothesis on an odd prime $p$:   

\vspace{0.1in} 
 
\noindent \textit{Hypothesis (*):} 
$p \nmid N$ and the extension $\Q(E[p])/\Q$ has Galois group isomorphic to $\GL_2(\F_p)$, i.e. the $\mod p$ Galois 
representation $\rho_{E, p} : \Gal(\Qbar / \Q) \ra \GL(E[p])$ is surjective.  

\vspace{0.1in}

\noindent For such a prime\footnote{Recently, Byungchul Cha (see~\cite{cha:vanishing}) has been able to 
weaken Hypothesis (*).} $p$, Kolyvagin gives a precise formula for the order 
of the $p$-primary part of $\Sha(E/K)$ by constructing explicit elements in the $p$-power Selmer group $\Sel_{p^\infty}(E/K)$ 
from Heegner points, namely 
\begin{equation*}
\# \Sha(E/K)[p^\infty] = p^{2(m_0 - m_\infty)}, \leqno{(1)}
\end{equation*}
where $m_0 = [E(K) : \Z y_K]$ and $m_\infty$ is a non-negative integer which is defined in terms of global 
$p$-divisibility of the Heegner points used for the construction of the classes. We give a precise 
definition of $m_\infty$ in Section~\ref{subsec:basic-heeg}. The above formula provides a very strong evidence for 
Conjecture~\ref{conj:GZ-BSD}. Yet, Kolyvagin's arguments give no indication of how to relate the correction factor 
$m_\infty$ to the Manin constant and the Tamagawa numbers. 

\begin{remark}
Since $p \nmid N$, a result of Abbess and Ullmo (see~\cite{abbess-ullmo:manin}) shows that $p \nmid c$. Thus, 
$\ord_p \left ( c \cdot \prod_{q \mid N} c_q \right ) = \ord_p \left ( \prod_{q \mid N} c_q \right )$ (see 
also~\cite{agashe-ribet-stein:manin-constant} for an account on the known results about $c$). 
\end{remark} 

Our main goal is to attempt to relate $m_\infty$ to the Tamagawa numbers of $E$. More precisely, by combining (1) with 
Conjecture~\ref{conj:GZ-BSD} and the above remark, we obtain the following reformulation of the Birch and Swinnerton-Dyer formula 
for $E_{/K}$

\begin{conjecture}[BSD conjectural formula]\label{conj:m_eq_tam}
If $p$ satisfies \emph{Hypothesis (*)} then 
$$
m_\infty = \ord_p \left ( \prod_{q \mid N} c_q \right ).
$$
\end{conjecture}

\noindent The main contribution of this paper is the following theorem which provides evidence for 
Conjecture~\ref{conj:m_eq_tam}: 

\begin{theorem}\label{thm:main}
Assume that $p$ satisfies \emph{Hypothesis (*)}. If $\ds m_{\max} = \max_{q \mid N} \ord_p(c_q)$, 
then $m_\infty \geq m_{\max}$.  
\end{theorem}

One obtains as an immediate consequence of the above theorem and Kolyvagin's formula for the order of $\Sha(E/K)[p^\infty]$ 
the following

\begin{corollary}\label{cor:sha}
The $p$-primary part of the Shafarevich-Tate $\Sha(E/K)$ satisfies 
$$
\#\Sha(E/K)[p^\infty] \leq p^{2m_0 - 2m_{\max}}, 
$$
where $\ds m_{\max} = \max_{q \mid N} \ord_p(c_q) $. In particular, if $p$ divides at most one Tamagawa number, the 
above upper bound coincides with the exact upper bound for $\# \Sha(E/K)[p^\infty]$ predicted by the Birch and Swinnerton-Dyer 
conjectural formula for $E_{/K}$.  
\end{corollary}

\begin{remark}
It is explained in~\cite[Thm.8]{stein-wuthrich} that one could use results of Schneider, Perrin-Riou and Kato 
from Iwasawa theory to obtain the exact upper bound on the $p$-primary part of the Shafarevich-Tate group for $E_{/\Q}$ predicted 
by the Birch and Swinnerton-Dyer conjecture, but such a bound would be conditional on the hypothesis that the $p$-adic height 
$\hat{h}_p(P) \ne 0$, where $P$ is a generator of the Mordell-Weil group\footnote{Here, $\hat{h}_p(P)$ is also known as the $p$-adic regulator.}. 
One suspects (see~\cite[Rem.0.13(ii)]{colmez:p-adic_bsd}) that the last hypothesis might be at least as difficult as proving Leopold's 
conjecture. Our Corollary~\ref{cor:sha} in the analogous situation over $K$ does not depend on such a hypothesis.      
\end{remark}

We prove the theorem by refining Kolyvagin's original arguments which he used to study the exact structure of Selmer groups. We also apply 
several techniques from the theory of Kolyvagin systems as developed by Mazur, Rubin and Howard (see~\cite{mazur-rubin:kolyvagin_systems} 
and~\cite{howard:heeg}). The paper is organized as follows: Section~\ref{sec:notation} contains some notations. In Section~\ref{sec:heegpts} 
we review the basics of Heegner points over ring class fields, define the numbers $m_0$ and $m_\infty$ which appear in Kolyvagin's 
formula for the order of $\Sha(E/K)[p^\infty]$ and recall the construction of Kolyvagin's classes. We also state some known results 
about the reduction of the Heegner points at bad places which are necessary for understanding the precise local behavior of the classes. 
Section~\ref{sec:selmod} is about Selmer modules as introduced by Mazur and Rubin. We discuss various local conditions to be used in the argument, 
recall the notion of a Selmer structure and the associated Selmer modules and introduce a new Selmer structures which will replace  
the standard Selmer structure obtained from the Kummer local conditions in Kolyvagin's arguments. The main contribution in this section is 
Proposition~\ref{prop:refine} where we obtain more refined local properties of the constructed classes than the ones implicitly used by Kolyvagin. 
Section~\ref{sec:global-duality} contains the necessary global duality results (consequences of the duality theorem of Poitou and Tate) on Galois 
cohomology. Finally, we prove our main theorem in Section~\ref{sec:main-thm}. The proof uses extensively a combination of \v{C}ebotarev density 
theorem with the global duality results. We first prove the result in an easy case (Theorem~\ref{thm:min-core}) and then use some techniques of 
Kolyvagin to reduce the general case to the easy case in Section~\ref{subsec:gencase}.        

\section*{Acknowledgements} 

I am thankful to William Stein for suggesting the problem to me and for the various helpful discussions. I am very grateful to my advisor Ken Ribet 
for his constant encouragement, guidance, support and numerous helpful conversations on the project.  

I am indebted to Christophe Cornut for his enormous help and support while working on this project, for the careful 
reading of the preliminary draft and for sharing his deep and fine understanding of Kolyvagin's method.    

I am also greatful to many people, especially Byungchul Cha, Mladen Dimitrov, Olivier Fouquet, 
Ralph Greenberg, Ben Howard, Grigor Grigorov, Barry Mazur, Robert Pollack, Bjorn Poonen, Karl Rubin and  
Olivier Wittenberg, for the many helpful conversations. It is also a pleasure to thank the Mathematical Sciences Research 
Institute at Berkeley and the Institut de Math\'{e}matiques de Jussieu, Paris for their kind hospitality while working on the 
project.

\section{Notations}\label{sec:notation}
Throughout the whole paper $\Kbar$ will denote a fixed algebraic closure of $K$, $\cO_K$ - the ring of integers of $K$, $\cO_c$ - the 
non-maximal order of conductor $c$ in $\cO_K$ (i.e., $\cO_c = \Z + c \cO_K$), $K[c]$ - the ring class field extension of $K$ of conductor $c$ (i.e., 
an abelian extension of $K$ whose Galois group is isomorphic to $\Pic(\cO_c)$ and which is Galois and dihedral over $\Q$; e.g., $K[1]$ is the Hilbert class 
field of $K$). For a number field $L$ 
and a place $w$ of $L$, $L_w$ will be the completion of $L$ at $w$, $\cO_w$ - the ring of integers of $L_w$, $\Lbar_w$ - a fixed algebraic closure of 
$L_w$ and $L_w^{\ur}$ - the maximal unramified extension of $L_w$. Whenever $F$ is a field (local or global) with a fixed algebraic closure $\overline{F}$, $G_F$ will denote the Galois group $\Gal(\overline{F} / F)$. 
If $M$ is a $G_F$-module with a continuous action of $G_F$, $\H^1(F, M)$ will be the Galois cohomology group $\H^1(G_F, M)$. 

Moreover, fix an embedding $\iota_v : \Kbar \hra \overline{K_v}$ for each place $v$ of $K$ (this corresponds to fixing a place of $\Kbar$ above 
$v$ for every place $v$ of $K$).  

We also accept Kolyvagin's original notation for abelian $p$-groups. If   
$A \isom \Z / p^{x_1}\Z \oplus \Z / p^{x_2} \Z \oplus \dots \oplus \Z / p^{x_n} \Z$  
with $x_1 \geq x_2 \geq \dots \geq x_n$, the sequence of invariants of $A$ is defined as 
$$
\Inv(A) = (x_1, x_2, \dots, x_n). 
$$
Also, for $a \in A$, $\ord (a)$ stands for the integer $n$, such that $p^n$ is the order of $a$ in $A$. 

\section{Heegner Points Over Ring Class Fields}\label{sec:heegpts}

Kolyvagin used Heegner points over ring class fields for certain non-maximal orders of $K$ to construct explicit 
cohomology classes in $\H^1(K, E[p^m])$ for each $m$ (see~\cite{kolyvagin:euler_systems},
~\cite{gross:kolyvagin} or~\cite{mccallum:kolyvagin}). He used the classes to study the structure of the   
Selmer group $\Sel_{p^\infty}(E/K)$ (see~\cite{kolyvagin:structure_of_selmer}). 

\subsection{Defining Heegner points over ring class fields}\label{subsec:basic-heeg}

\vspace{0.1in}

\noindent \textit{1. Heegner points over ring class fields.} 
Let $\cO_c = \Z + c\cO_K$ be the order of conductor $c$ in $\cO_K$ and let $\cN_c = \cN \cap \cO_c$. 
Then $\cN_c \subset \cO_c$ is an invertible ideal, $\cO_c / \cN_c \simeq \Z/N\Z$ and the map $\C / \cO_c \ra \C / \cN_c^{-1}$ is 
a cyclic isogeny of degree $N$, so it defines a point $x_c \in X_0(N)(\C)$ which is $K[c]$-rational by the theory of complex multiplication.  
One can use the parameterization $\varphi : X_0(N) \ra E$ to construct a point $y_c = \varphi (x_c) \in E(K[c])$. 

\vspace{0.1in}

\noindent \textit{2. Kolyvagin primes and conductors.} We call a prime number $\ell$ a \emph{Kolyvagin prime} 
relative to $E$, $K$ and $p$ if $\ell$ is inert in $K$ and $p$ divides both $a_\ell$ and $\ell+1$. 
For each Kolyvagin prime $\ell$ let $M(\ell) = \ord_p(a_\ell, \ell+1)$.  Denote by $\Lambda^r$ the set 
of all square-free products of exactly $r$ Kolyvagin primes and by $\ds \Lambda = \bigcup_r \Lambda^r$ 
(by convention, $\Lambda^0 = \{1\}$). For each $c \in \Lambda$ define $\ds M(c) = \min_{\ell \mid c} M(\ell)$. 
For the purpose of our argument, we will also need to consider the subset $\Lambda_m^r \subset \Lambda^r$ defined 
as  
$$
\Lambda_m^r = \{c \in \Lambda^r : M(c) \geq m\}. 
$$
We also set $\ds \Lambda_m = \bigcup_r \Lambda_m^r$. 

\vspace{0.1in}

\noindent \textit{3. Heegner points relations.} The Heegner points satisfy two important relations 
(see~\cite[Prop.3.7]{gross:kolyvagin} for the derivations). 

\begin{enumerate}
\item \emph{Distribution relation.} For every $c \in \Lambda$ and $\ell \in \Lambda^1$, 
$$
\Tr_{K[c\ell]/K[c]} y_{c\ell} = a_\ell y_c. 
$$

\item \emph{Congruence relation.} For any prime $\lambda_{c\ell}$ of $K[c\ell]$ above $\lambda \mid \ell$, 
$$
y_{c\ell} \equiv \Fr_{\lambda_c} y_c \mod \lambda_{c\ell}, 
$$
where $\lambda_c$ is the prime of $K[c]$ below $\lambda_{c\ell}$ and $\Fr_{\lambda_c}$ is the associated 
Frobenius. 
\end{enumerate}

\vspace{0.1in}

\noindent \textit{4. Kolyvagin derivative operators.} 
Let $\cG_c= \Gal(K[c]/K)$ and $G_c = \Gal(K[c]/K[1])$. 
For each $\ell \in \Lambda^1$, the group $G_\ell$ is cyclic of order 
$\ell + 1$. Indeed, 
$$
G_\ell \simeq (\cO_K / \ell \cO_K)^\times / (\Z / \ell \Z)^\times \simeq \F_\lambda^\times / \F_\ell^\times. 
$$
Moreover, $\ds G_c \isom \prod_{\ell \mid c} G_\ell$ (to see this, one uses that the subgroup of $G_c$ fixing $K[c/\ell]$ is isomorphic to $G_\ell$). 
Next, fix a generator $\sigma_\ell$ of $G_\ell$ for each $\ell \in \Lambda^1$. Define 
$D_\ell = \sum_{i = 1}^\ell i \cdot \sigma_\ell^i \in \Z[G_\ell]$ and let 
$$
D_c = \prod_{\ell \mid c} D_\ell \in \Z[G_c]. 
$$
Note that $(\sigma_\ell - 1) D_\ell = 1 + \ell - \Tr_{K[\ell]/K[1]}$. 

We refer to $D_c$ as the \emph{Kolyvagin derivative operators}. Finally, let $S$ be a set of coset representatives 
for the subgroup $G_c \subseteq \cG_c$. Define 
$$
P_c = \sum_{s \in S} s D_c y_c \in E(K[c]).  
$$ 
The points $P_c$ are derived from the points $y_c$, so we will refer to them as \emph{derived Heegner points}. 

\vspace{0.1in}

\noindent \textit{5. Defining $m_0$ and $m_\infty$.} 
Let $m'(c)$ be the largest 
positive integer, such that $P_c \in p^{m'(c)}E(K[c])$ (if $P_c$ is torsion then $m'(c) = \infty$). Define a 
function $m : \Lambda \ra \Z$ by 
$$
m(c) = \left \{ 
\begin{array}{ll} 
m'(c) & \textrm{ if }m'(c) \leq M(c) \\
\infty & \textrm{ else }
\end{array} \right . 
$$ 
Finally, let $\ds m_r = \min_{c \in \Lambda^r} m(c)$. Note that by our convention for $\Lambda^0$ and by the fact that $P_1 = y_K$, 
$m_0 = \ord_p[E(K) : \Z y_K] < \infty$ for $y_K$ has infinite order and $E(K)$ has rank one by the result of Kolyvagin. 

Kolyvagin considered the sequence $\{m_r\}$ and proved (see~\cite[Thm.C]{kolyvagin:structureofsha}) that $m_r \geq m_{r+1}$ 
for every $r \geq 0$, hence, it stabilizes to a limit $m_\infty$, which, in our particular setting, is always finite since $m_0$ is 
finite\footnote{In fact, Kolyvagin proved that $m_r \geq m_{r+1}$ without the assumption that the Heegner point $y_K = P_1$ has infinite order in 
$E(K)$. In this situation, one could still define $m_\infty$, but it is not at all obvious whether $m_\infty < \infty$. Kolyvagin conjectured this 
for all elliptic curves (see~\cite[Conj.C]{kolyvagin:structure_of_selmer} for the original statement of the conjecture 
and~\cite{jetchev-lauter-stein:kolconj} for some applications of Kolyvagin's conjecture and some computational and theoretical evidence.}.   

\vspace{0.1in}

\noindent \textit{6. Kolyvagin classes $\kappa_{c, m} \in \H^1(K, E[p^m])$.} Let $c \in \Lambda_m$. To construct the class 
$\kappa_{c, m} \in \H^1(\cF, E[p^m])$ one first observes that the image $\tilde{P}_c$ of $P_c$ in $E(K[c])/p^m E(K[c])$ is 
fixed by $\cG_c$ (see~\cite[Prop.3.6]{gross:kolyvagin}). Since the Galois representation $\rho_{E, p}$ is surjective, 
the restriction map $\H^1(K, E[p^m]) \xra{\res} \H^1(K[c], E[p^{m}])^{\cG_c}$ is an isomorphism (see~\cite[pp.241-242]{gross:kolyvagin}), 
so if $\delta_c : E(K[c])/p^{m} E(K[c]) \ra \H^1(K[c], E[p^m])$ 
is the Kummer map, one can define
$$
\kappa_{c,m} := \res^{-1}(\delta_c(\tilde{P}_c)) \in \H^1(K, E[p^m]).  
$$ 
It follows from the definition of the Kummer map that $\kappa_{c, m} = 0$ if and only if $P_c \in p^m E(K[c])$ (which is equivalent to 
$m \leq m(c)$). Moreover, if $m > m(c)$ then $\ord (\kappa_{c, m}) = m - m(c)$. 
The class $\kappa_{c, m}$ is represented by the 1-cocycle
\begin{eqnarray}\label{eqn:cocycle}
\sigma \mapsto \sigma \left (\frac{P_c}{p^m} \right ) - \frac{P_c}{p^m} - \frac{(\sigma - 1)P_c}{p^m}, 
\end{eqnarray}
where $\ds \frac{(\sigma - 1)P_c}{p^m}$ is the unique $p^m$-division point of $(\sigma-1)P_c$ in $E(K[c])$ (see~\cite[Lem.~4.1]{mccallum:kolyvagin}). 

Finally, let $\eps = \pm 1$ be the eigenvalue of the Atkin-Lehner (Fricke) involution $w_N$ on the eigenform $f$ 
corresponding to $E$, i.e. $f | w_N = \eps \cdot f$. For each $c \in \Lambda_m$, let 
$\eps(c) = \eps \cdot (-1)^{f_c}$ where $f_c = \# \{\ell : \ell \mid c\}$. It follows from~\cite[Prop.5.4(ii)]{gross:kolyvagin} that 
$\kappa_{c, m}$ lies in the $\eps(c)$-eigenspace for the action of complex conjugation on $\H^1(K, E[p^m])$.

\vspace{0.1in} 

\noindent \textit{7. Kolyvagin classes $\tilde{\kappa}_{c, m} \in \H^1(K, E[p^m])$.} Suppose that $c \in \Lambda$ satisfies 
$m + m(c) \leq M(c)$. We construct a class $\tilde{\kappa}_{c, m} \in \H^1(K, E[p^m])$, such that the cyclic 
$\Z / p^m \Z$-submodule generated by $\tilde{\kappa}_{c, m}$ is free of rank one and contains the original class $\kappa_{c, m}$.  
Indeed, consider the short-exact sequence 
$$
0 \ra E[p^m] \ra E[p^{m+m(c)}] \xra{p^{m}} E[p^{m(c)}] \ra 0. 
$$
The corresponding long-exact sequence on Galois cohomology gives an exact sequence  
$$
0 \ra \H^1(K, E[p^{m}]) \hra \H^1(K, E[p^{m+m(c)}]) \xra{p^{m}} \H^1(K, E[p^{m(c)}]),  
$$
since $\H^0(K, E[p^{m(c)}]) = E(K)[p^{m(c)}] = 0$. Since $m + m(c) \leq M(c)$, one can consider the cohomology class 
$\kappa_{c, m + m(c)} \in \H^1(K, E[p^{m+m(c)}])$. It follows from the definition of  
$\kappa_{c, m}$ that $p^{m(c)} \kappa_{c, m + m(c)} = \kappa_{c, m}$ in $\H^1(K, E[p^{m+m(c)}])$. We claim that 
$\kappa_{c, m+m(c)}$ is in the image of $\H^1(K, E[p^m])$ under the above inclusion. Indeed, since $\ord( \kappa_{c, m+m(c)} ) = m$ then it is 
in the kernel of the second map in the above exact sequence, i.e., it comes from a class $\tilde{\kappa}_{c, m} \in \H^1(K, E[p^m])$. Moreover, 
$$
\ord \tilde{\kappa}_{c, m}= \ord \kappa_{c, m + m(c)} = m,
$$ 
i.e., $\tilde{\kappa}_{c, m}$ spans a free $\Z / p^m \Z$-submodule of $\H^1(K, E[p^m])$.

\subsection{Reduction properties of Heegner points}\label{subsec:heegred}
Let $J = \Jac(X_0(N))$ be the Jacobian of the modular curve $X_0(N)$. The following proposition summarizes several results 
discussed in the paper of Gross and Zagier (for details, see~\cite[\S III.1]{gross-zagier} and~\cite[\S III.3,~Prop.3.1]{gross-zagier}).   

\begin{proposition}\label{prop:heegred}
Let $v$ be a prime in $K$ which is a prime of bad reduction for $E$, i.e. $v \mid N$. For any 
conductor $c$ which is prime to $\chr(v)$ and any prime $w \mid v$ of $K[c]$, the Heegner divisor class 
$[(x_c) - (\infty)] \in J(K[c]_w)$ lies, up to translation by the rational divisor class 
$[(0) - (\infty)]$, in $J^0(K[c]_w)$, where $J^0(K[c]_w) \subseteq J(K[c]_w)$ is the subgroup of all points which specialize to the 
identity component of the N\'eron model of $J$.   
\end{proposition}

Since the image of the divisor class $[(0) - (\infty)]$ under $\vphi$ is a rational torsion point on $E$, we obtain the 
following 

\begin{corollary}\label{cor:red}
The Heegner point $y_c$ lies, up to translation by a rational torsion point of $E$, on $E^0(K[c]_w)$, where 
$E^0(K[c]_w))$ is the subgroup of $E(K[c]_w)$ of the points which specialize to the identity component of 
the N\'eron model of $E$.  
\end{corollary}

\begin{remark}
Kolyvagin used these reduction properties show that the classes that he constructed satisfy certain local conditions 
at $v$ (see~\cite[Prop.6.2]{gross:kolyvagin} for the application). We will use the above corollary in the next section to 
prove that Kolyvagin's classes satisfy certain refined local conditions which will be a key observation towards proving 
Theorem~\ref{thm:main}.    
\end{remark}

\section{Selmer Modules}\label{sec:selmod}

\subsection{Local Selmer conditions}
Fix a positive integer $m$ and a place $v$ of $K$. We define and compare several local Selmer conditions. 

\vspace{0.1in}

\noindent \emph{1. Unramified local Selmer condition.} The unramified local Selmer condition for $E[p^m]$ is defined to be the 
subgroup  
$$
\H^1_{\ur}(K_v, E[p^m]) := \ker \{ \H^1(K_v, E[p^m]) \ra \H^1(K_v^{\ur}, E[p^m])\}. 
$$

\begin{remark}\label{rem:fin-ur}
If the module $E[p^m]$ is unramified at $v$ and $v \nmid p$ then one also refers to the unramified local Selmer 
condition as the \emph{finite local condition}\footnote{We use the alternative terminology to stick to the 
standard terminology introduced in~\cite{rubin:book} and~\cite{mazur-rubin:kolyvagin_systems}. For the general definition of the 
finite condition on the Tate module $T$, on $V = T \otimes \Q_p$ and on the finite quotients of $T$ see~\cite[\S 1.3]{rubin:book}. It 
is known that if $T$ is unramified then $\H^1_{\f}(K_v, T) = \H^1_{\ur}(K_v, T)$, but in general, $\H^1_{\ur}(K_v, T) \subseteq \H^1_{\f}(K_v, T)$ 
with finite index. One should note that the finite local condition is used in the statement of the Bloch-Kato conjectures (see~\cite{bloch-kato}).} and 
denote it by $\H^1_{\f}(K_v, E[p^m])$. 

\comment{We will often use the \emph{singular quotient} 
$\H^1_{\s}(K_v, E[p^m]) = \H^1(K_v, E[p^m])/\H^1_{\f}(K_v, E[p^m])$ together with the exact 
sequence 
$$
0 \ra \H^1_{\f}(K_v, E[p^m]) \ra \H^1(K_v, E[p^m]) \ra \H^1_{\s}(K_v, E[p^m]) \ra 0.   
$$  
}
\end{remark}

\vspace{0.1in}

\noindent \emph{2. Transverse local Selmer condition.} For a Kolyvagin prime $\ell \in \Lambda^1$ we denote by $K[\ell]_\lambda$ the completion of  
$K[\ell]$ at the place below the place corresponding to the fixed embedding $\iota_\lambda : \Kbar \hra \overline{K_\lambda}$. 
Define the transverse condition as 
$$
\H^1_{\tr}(K_v, E[p^m]) := \ker \{ \H^1(K_v, E[p^m]) \ra \H^1(K[\ell]_\lambda, E[p^m])\}. 
$$
We give a more explicit description of the transverse condition in the next section of the paper. 
 
\comment{
\noindent \emph{2. $L$-transverse local Selmer condition.} For a maximal totally and tamely ramified abelian extension of $K_v$ 
(e.g., $L = K[\ell]_{\lambda'}$ for $\ell \in \Lambda^1_m$ and $\lambda'$ - any prime of $K[\ell]$ above $\ell$), define the $L$-transverse 
condition as 
$$
\H^1_{L-\tr}(K_v, E[p^m]) := \ker \{ \H^1(K_v, E[p^m]) \ra \H^1(L, E[p^m])\}. 
$$
}

\vspace{0.1in}

\noindent \emph{3. Kummer local condition.} The local condition $\H^1_{\Kum}(K_v, E[p^m])$ is defined to be the image of the local 
Kummer map 
$$
\delta_v : E(K_v) / p^m E(K_v) \hra \H^1(K_v, E[p^m]) 
$$ 
for every non-archimedian place $v$ of $K$. 
\comment{The local condition at the complex archimedian place is defined to be $\H^1(\C, E[p^m])$.}

\vspace{0.1in} 

\noindent \emph{4. Stringent Kummer local condition.} We denote it by $\H^1_{\Kum^0}(K_v, E[p^m])$ and define it as the image of $E^0(K_v)$ under 
the map      
$$
E(K_v) \twoheadrightarrow E(K_v)/p^m E(K_v) \xra{\delta_v} \H^1(K_v, E[p^m]). 
$$  

\vspace{0.1in}

The next proposition compares the unramified local condition with the Kummer local condition in the 
case when $v$ is a place of good reduction for $E$ and $v \nmid p$.  
 
\begin{proposition}\label{prop:ur-pts}
Let $v$ be a place of good reduction for $E$, such that $v \nmid p$. Then  
$\H^1_{\ur}(K_v, E[p^m]) = \H^1_{\Kum}(K_v, E[p^m])$.  
\end{proposition}

We need the following two lemmas for the proof of the above propositions. 

\begin{lemma}\label{lemma:loc-ur-pts}
Let $v \nmid p$ is any finite place of $K$ then the group $E^0(K_v^{\ur})$ is $p$-divisible. 
\end{lemma}

\begin{proof}
Consider the diagram  
$$
\xymatrix{
0 \ar[r] & E^1(K_v^{\ur}) \ar[r]\ar[d]^{\cdot p}& E^0(K_v^{\ur}) \ar[r]\ar[d]^{\cdot p} & \tilde{E}_{ns}(\kbar_v) \ar[r]\ar[d]^{\cdot p} & 0 \\
0 \ar[r] & E^1(K_v^{\ur}) \ar[r] & E^0(K_v^{\ur}) \ar[r] & \tilde{E}_{ns}(\kbar_v) \ar[r] & 0.  
} 
$$ 
Let $L_w / K_v$ be a finite unramified extension of $K_v$ ($w$ is a place above $v$). Then $E^1(L_w) \isom \widehat{E}(\m_w)$ (see~\cite[Ch.VII,~Prop.2.2]{silverman:aec}), where $\widehat{E}$ is the formal group associated to $E$ and $\m_w$ is the maximal ideal 
of the ring of integers $\cO_w$ of $L_w$. Since $\chr(v) \ne p$, multiplication-by-$p$ is an isomorphism on the formal group 
according to~\cite[Ch.IV,\S3]{silverman:aec}, i.e. the left vertical map is an isomorphism by taking a limit 
over all finite unramified extensions of $K_v$. The right vertical map is surjective since $\kbar_v$ is 
algebraically closed and $\tilde{E}_{ns}$ is either an elliptic curve over $k_v$ or $\Gm$, or $\mathbb{G}_a$. 
Thus, the map in the middle is surjective by the snake lemma, i.e., $E^0(K_v^{\ur})$ is $p$-divisible.  
\end{proof}

\begin{remark}
The above short proof uses Weierstrass minimal models. There is another proof which uses   
N\'eron models (see~\cite[\S 7.3]{neronmodels}) and which works for abelian varieties of arbitrary 
dimensions. 
\end{remark}

\comment{
\begin{lemma}\label{lemma:h1triv}
If $v \nmid p$ is a place of good reduction for $E$ then 
$$
\H^1(K_v^{\ur}/K_v, E(K_v^{\ur})) = 0.  
$$
\end{lemma}
 
\begin{proof}
According to~\cite[I.3.8]{milne:duality}, 
$$
\H^1(K_v^{\ur}/K_v, E(K_v^{\ur})) \isom \H^1(\kbar_v / k_v, \Phi_{E, v}(\kbar_v)). 
$$ 
The last group is trivial since the component group $\Phi_{E,v}$ is trivial. 
\end{proof} 
} 
 
\begin{proof}[Proof of Proposition~\ref{prop:ur-pts}]
Consider the following commutative diagram 
$$
\xymatrix{
 & & \H^1(K_v^{\ur}/K_v, E(K_v^{\ur}))[p^m] \ar[d]& & \\
E(K_v) / p^m E(K_v) \ar[r] & \H^1(K_v, E[p^m]) \ar[r]\ar[d]\ar@{-->}[dr]^{\theta}& \H^1(K_v, E)[p^m] \ar[d]\ar[r]& 0 \\
E(K_v^{\ur}) / p^m E(K_v^{\ur})\ar[r] & \H^1(K_v^{\ur}, E[p^m]) \ar[r]& \H^1(K_v^{\ur}, E)[p^m] \ar[r]& 0     
} 
$$
Lemma~\ref{lemma:loc-ur-pts} and the Kummer sequence for $E/K_v^{\ur}$ imply that 
$\H^1_{\ur}(K_v, E[p^m]) = \ker ({\theta})$. 

Similarly, Lang's theorem (see~\cite{lang:finitefields}) implies that $\H^1(K_v^{\ur}/K_v, E(K_v^{\ur})) = 0$, so it follows from the Kummer 
sequence for $E/K_v$ that the image of $E(K_v)/p^m E(K_v)$ under the Kummer map is equal to $\ker (\theta)$. This proves the proposition. 
\end{proof}

\subsection{Local conditions at Kolyvagin primes}
Let $\ell \in \Lambda_m^1$ be a Kolyvagin prime and $\lambda \mid \ell$ be the unique prime of $K$ above $\lambda$. The main property of these 
primes which we will use in this section is the fact that $G_{K_\lambda}$ acts trivially on $E[p^m]$. Indeed, $G_{K_\lambda}$ acts through its 
quotient $\Gal(K_\lambda^{\ur}/K_\lambda)$ which is topologically generated by $\Fr_\lambda$. But $\Fr_\lambda = \Fr_\ell^2$ and since $\ell \in \Lambda^1_m$ 
then $\Fr_\ell$ acts as complex conjugation on $E[p^m]$, i.e., $\Fr_\lambda$ acts trivially.  

\vspace{0.1in}

\noindent \textit{1. Local cohomology groups.} We have 
$$
\H^1(K_\lambda, E[p^m]) = \Hom(G_{K_\lambda}, E[p^m]) \isom \Hom(G_{K_\lambda}^{\ab}/p^m, E[p^m]). 
$$
By local class field theory, 
$$
G_{K_\lambda}^{\ab}/p^m \isom \Gal(K_\lambda^{\ur}/K_\lambda)/p^m \oplus \Gal(K[\ell]_\lambda / K_\lambda)/p^m. 
$$
Indeed, 
$$
G_{K_\lambda}^{\ab}/p^m \isom K_\lambda^{\times} / (K_\lambda^\times)^{p^m} = \ell^{\Z /p^m\Z} \oplus U_\lambda / (U_\lambda)^{p^m}, 
$$
where $U_\lambda \subset \cO_\lambda$ is the group of units. Local class field theory then gives us isomorphisms 
$$
\ell^{\Z /p^m\Z} \isom \Gal(K_\lambda^{\ur}/K_\lambda)/p^m
$$
and 
$$
U_\lambda / (U_\lambda)^{p^m} \isom \Gal(K[\ell]_\lambda / K_\lambda)/p^m, 
$$
so we obtain the above decomposition. In particular, 
$$
\H^1(K_\lambda, E[p^m]) = \H^1_{\ur}(K_\lambda, E[p^m]) \oplus \H^1_{\tr}(K_\lambda, E[p^m]), 
$$
where 
$$
\H^1_{\ur}(K_\lambda, E[p^m]) = \{f \in \Hom(G_{K_\lambda}^{\ab}/p^m, E[p^m]) : f(\Gal(K[\ell]_\lambda / K_\lambda)/p^m) = 0\}
$$
and 
$$
\H^1_{\tr}(K_\lambda, E[p^m]) = \{f \in \Hom(G_{K_\lambda}^{\ab}/p^m, E[p^m]) : f(\Gal(K_\lambda^{\ur} / K_\lambda)/p^m) = 0\}. 
$$
A homomorphism in $\H^1_{\ur}(K_\lambda, E[p^m])$ is determined by the image of $\Fr_\lambda$, whereas a homomorphism in 
$\H^1_{\tr}(K_\lambda, E[p^m])$ is determined by the image of the restriction to the decomposition group of the fixed generator $\sigma_\ell$ 
of the cyclic group $\Gal(K[\ell]/K[1])$. Thus,  
$$
\H^1_{\ur}(K_\lambda, E[p^m]) \isom E[p^m] \textrm{ and } \H^1_{\tr}(K_\lambda, E[p^m]) \isom E[p^m]. 
$$
Since $\tau \in \Gal(K_\lambda/\Q_\ell)$ acts on $\Hom(G_{K_\lambda}^{\ab}, E[p^m])$ by sending $f \mapsto f^\tau$, where 
$f^{\tau}(\sigma) = f(\tau^{-1}\sigma \tau)^{\tau}$ for $\sigma \in G_{K_\lambda}^{\ab}$ then the first of the above 
isomorphism is Galois equivariant and the second is anti-equivariant (since $K[\ell] / \Q$ is dihedral). Therefore, 
$$
\H^1_{\ur}(K_\lambda, E[p^m])^{\pm} \isom E[p^m]^{\pm} \textrm{ and } \H^1_{\tr}(K_\lambda, E[p^m])^{\pm} \isom E[p^m]^{\mp}. 
$$
The existence of the Weil pairing implies that $E[p^m]$ splits into two eigenspaces of complex conjugation each of which is free of rank 
one over $\Z/p^m \Z$ and thus, $\H^1_{\ur}(K_\lambda, E[p^m])^{\pm}$ and $\H^1_{\tr}(K_\lambda, E[p^m])^{\pm}$ are all free of rank one over 
$\Z/p^m \Z$. 

Finally, one can use the explicit description of the $\pm$-eigenspaces to conclude that $\H^1_{\tr}(K_\lambda, E[p^m])$ is self-dual with respect to the 
Tate local pairing\footnote{For a different proof using class field theory and Kummer theory see~\cite[Prop.1.3.2(ii)]{mazur-rubin:kolyvagin_systems}.}. Indeed, it suffices to show that $\H^1_{\tr}(K_\lambda, E[p^m])$ is a self-orthogonal (since it has the dimension of a maximal isotropic 
subspace). Since $\H^1_{\tr}(K_\lambda, E[p^m])^{\pm}$ are both cyclic $\Z/p^m\Z$-modules, each of them is self-orthogonal, so it will be enough to show that 
$\H^1_{\tr}(K_\lambda, E[p^m])^{+}$ is orthogonal to $\H^1_{\tr}(K_\lambda, E[p^m])^{-}$. The last is an immediate consequence of the 
$\Gal(K_\lambda /\Q_\ell)$-equivariancy of the Tate local pairing (i.e., $\langle \tau x,  \tau y\rangle_\lambda = \langle x, y\rangle_\lambda$). The 
same argument shows that $\H^1_{\ur}(K_\lambda, E[p^m])$ is self-dual as well.

\vspace{0.1in}

\noindent \textit{2. A comparison homomorphism.}  
We constructs a comparison homomorphism (non-canonically) between the unramified and the transverse local conditions         
$$
\phi_\lambda : \H^1_{\ur}(K_\lambda, E[p^m]) \ra \H^1_{\tr}(K_\lambda, E[p^m]).    
$$
Here, anti-equivariant means that  
$$
\phi_\lambda \left (\H^1_{\ur}(K_\lambda, E[p^m])^{\pm} \right ) = \H^1_{\tr}(K_\lambda, E[p^m])^{\mp}.
$$ 
We saw in 2. that $\H^1_{\ur}(K_\lambda, E[p^m]) \isom E[p^m]$ and in 3. we obtained a non-canonical identification of  
$\H^1_{\tr}(K_\lambda, E[p^m])$ with $E[p^m]$.     

We next construct an automorphism $\chi_\ell : E[p^{M(\ell)}] \ra E[p^{M(\ell)}]$ (see~\cite[\S 1.7]{howard:heeg}). First,  
$\Fr_\ell$ splits $\tilde{E}(k_\lambda)[p^\infty]$ into two eigenspaces which are cyclic groups of orders $p^{\ord_p( \ell + 1 - a_\ell )}$ 
and $p^{\ord_p(\ell + 1 + a_\ell)}$. Next, consider the following sequence of homomorphisms   
$$
E(K_\lambda) \ra \tilde{E}(k_\lambda) \ra \tilde{E}(k_\lambda)[p^\infty] \xra{p^{-M(\ell)}(a_\ell - (\ell+1)\Fr_\ell)} \tilde{E}(k_\lambda)[p^{M(\ell)}] \xra{\sim} E(K_\lambda)[p^{M(\ell)}],   
$$  
where the last map is the canonical lift. The kernel of this map (again, by the eigenspace decomposition) is $p^{M_{\ell}} E(K_\lambda)$. Thus, we obtain 
an isomorphism 
$$
\chi'_\ell : E(K_\lambda) / p^{M(\ell)} E(K_\lambda) \xra{\sim} E(K_\lambda)[p^{M(\ell)}]. 
$$ 
But 
$$
E(K_\lambda) / p^{M(\ell)} E(K_\lambda) = \H^1_{\ur}(K_\lambda, E[p^{M(\ell)}]) = E[p^{M(\ell)}] / (\Fr_v - 1) E[p^{M(\ell)}] = E[p^{M(\ell)}], 
$$
since $\Fr_\lambda = \Fr_\ell^2$ acts trivially on $E[p^{M(\ell)}]$. Therefore, we obtain an isomorphism $\chi_\ell : E[p^{M(\ell)}] \ra E[p^{M(\ell)}]$, so 
one can define the comparison isomorphism as 
$$
\phi_\lambda : \H^1_{\ur}(K, E[p^m]) \isom E[p^m] \xra{\chi_\ell} E[p^m] \isom \H^1_{\tr}(K_\lambda, E[p^m]).  
$$
Since $\chi_\ell$ is equivariant and since $\H^1_{\tr}(K_\lambda, E[p^m])^{\pm} \isom E[p^m]^{\mp}$ then $\phi_\lambda$ is anti-equivariant. 

\begin{remark}
The construction resembles the finite-singular homomorphism construction of Mazur, Rubin (see~\cite[\S 1.2]{mazur-rubin:kolyvagin_systems}). 
The difference between our construction and the construction of Howard (see~\cite[Defn.1.1.8]{howard:heeg}) is that we have an extra twist by the automorphism $\chi_\ell$ for we avoid some technical difficulties when we establish the local relations at Kolyvagin primes between the explicit Kolyvagin classes. 
\end{remark}

\vspace{0.1in}

\noindent \emph{5. Comparison between $\kappa_{c,m}$ and $\kappa_{c\ell, m}$ at $\lambda$.} 
We will prove the following 

\begin{proposition}\label{prop:comparison}
For any $c \in \Lambda_m$ and $\ell \in \Lambda_m^1$ for which $\ell \nmid c$, 
$$
\phi_\lambda(\loc_\lambda(\kappa_{c, m})) = \loc_\lambda(\kappa_{c\ell, m}). 
$$
\end{proposition}

\begin{proof}
Using the definition of $\phi_\lambda$ in 4., we reduce the statement to showing that  
$$
\chi_\ell (\kappa_{c, m}(\Fr_\lambda)) = \kappa_{c\ell, m}(\sigma_\ell), 
$$
where (by abuse of notation) $\kappa_{c, m}$ and $\kappa_{c\ell, m}$ are the explicit cocycles 
(\ref{eqn:cocycle}). This equality will follow if we show that 
$$
\kappa_{c\ell, m} (\sigma_\ell) \equiv \frac{a_\ell-(1+\ell)\Fr_\ell}{p^m} P_c \mod \lambda', 
$$ 
where $\lambda'$ is the prime of $K[\ell]$ above $\ell$ discussed in 3. 

To verify the congruence, we first evaluate the cocycle (\ref{eqn:cocycle}) at $\sigma_\ell$
$$
\kappa_{c\ell, m}(\sigma_\ell) = \sigma_\ell \left ( \frac{P_{c\ell}}{p^m}\right ) - \frac{P_{c\ell}}{p^m} - 
\frac{(\sigma_\ell-1)P_{c\ell}}{p^m}. 
$$
Note that the value is defined over $K[\ell]_{\lambda'}^{\ur}$. Since $\sigma_\ell$ is an element of the inertia 
group at $\lambda$ (it generates a totally ramified extension), the point 
$\ds \sigma_\ell \left ( \frac{P_{c\ell}}{p^m} \right ) - \frac{P_{c\ell}}{p^m}$ reduces 
to zero $\mod \lambda'$ i.e.,  
$$
\kappa_{c\ell, m}(\sigma_\ell) \equiv -\frac{(\sigma_\ell - 1)P_{c\ell}}{p^m} \mod \lambda'. 
$$
Next, we use the Heegner points distribution and congruence relations (see Section~\ref{subsec:basic-heeg}) to obtain 
\begin{eqnarray*}
\frac{(\sigma_\ell - 1)P_{c\ell}}{p^m} &=& \frac{(\sigma_\ell - 1)\sum_{\sigma \in S}  \sigma D_\ell D_c y_{c\ell}}{p^m} = \frac{(\sigma_\ell - 1) D_\ell \sum_{\sigma \in S} \sigma D_c y_{c\ell}}{p^m} = \\
&=& \frac{(1+\ell - \Tr_{K[c\ell]/K[c]}) \sum_{\sigma \in S}  \sigma D_c y_{c\ell}}{p^m} = \\
&=& \sum_{\sigma \in S}  \sigma D_c  \left ( \frac{(1+\ell )y_{c\ell}}{p^m} - 
\frac{\Tr_{K[c\ell]/K[c]}y_{c\ell}}{p^m} \right ) = \\
&=& \sum_{\sigma \in S} \sigma D_c \left ( \frac{(1+\ell)y_{c\ell}}{p^m} - \frac{a_\ell y_c}{p^m}\right ) \equiv \\
&\equiv&  \frac{(1+\ell)\Fr_\ell - a_\ell}{p^m} \sum_{\sigma \in S} \sigma D_c  y_c = 
\frac{(1+\ell)\Fr_\ell - a_\ell}{p^m} P_c
 \mod \lambda'
\end{eqnarray*}
Thus, 
$$
\kappa_{c\ell, m} (\sigma_\ell) \equiv \frac{ a_\ell - (1+\ell)\Fr_\ell}{p^m} P_c \mod \lambda', 
$$
which completes the proof.
\end{proof}

\subsection{Selmer structures and Selmer modules}

Selmer structures and Selmer are discussed in great generality by Mazur and Rubin 
(see~\cite[Ch.2]{mazur-rubin:kolyvagin_systems}). Here, we only need to consider Selmer structures on the Galois modules 
$E[p^m]$ for various $m$.   

\vspace{0.1in}

\noindent \textit{1. Selmer structures.} A \emph{Selmer structure} $\cF$ on $E[p^m]$ consists of a choice of a local Selmer condition 
$\H^1_{\cF}(K_v, E[p^m]) \subseteq \H^1(K_v, E[p^m])$ for each place $v$ of $K$, such that for all, but finitely many $v$, 
$\H^1_{\cF}(K_v, E[p^m]) = \H^1_{\ur}(K_v, E[p^m])$.  

\vspace{0.1in}

\noindent \textit{2. Selmer modules.} Given a Selmer structure $\cF$ on $E[p^m]$, one can define the corresponding 
\emph{Selmer module} as 
$$
\H^1_{\cF}(K, E[p^m]) := \Ker \left \{\H^1(K, E[p^m]) \ra \bigoplus_{v}\H^1(K_v, E[p^m]) / 
\H^1_{\cF}(K_v, E[p^m]) \right \}, 
$$
where the sum is taken over all places $v$ of $K$. Since $\Gal(K/\Q)$ acts on $\H^1(K, E[p^m])$ ($E$ is defined over $\Q$), the Selmer module 
$\H^1_{\cF}(K, E[p^m])$ will be a stable $\Gal(K/\Q)$-subspace of $\H^1(K, E[p^m])$ provided $\ds \bigoplus_{v \mid q}\H^1_{\cF}(K_v, E[p^m])$ is a 
$\Gal(K/\Q)$-stable subspace of $\ds \bigoplus_{v \mid q} \H^1(K_v, E[p^m])$ for every rational prime $q$. This will always be the case in the 
sequel of this paper.  

\vspace{0.1in}

\noindent \textit{3. Modified Selmer structures.} Let $\cF$ be a Selmer structure on $E[p^m]$ and $a, b, c$ be integers, such 
that $abc \in \Lambda_m$. The \emph{modified Selmer structure $\cF^a_b(c)$} is the structure whose local conditions are obtained 
from those of $\cF$ by simply replacing them at the places $v \mid abc$ as follows:   
\begin{itemize}
\item If $v \mid c$ then $\H^1_{\cF^a_b(c)}(K_v, E[p^m]) = \H^1_{\tr}(K_v, E[p^m])$.  

\item If $v \mid a$ then $\H^1_{\cF^a_b(c)}(K_v, E[p^m]) = \H^1(K_v, E[p^m])$. 

\item If $v \mid b$ then $\H^1_{\cF^a_b(c)}(K_v, E[p^m]) = 0$.  
\end{itemize} 

\noindent Note that if $a = 1$ then we usually omit the $a$ and write $\cF_b(c)$. If $b = 1$, we omit the $b$ and write 
$\cF^a(c)$. 

\vspace{0.1in}

\noindent \textit{4. Dual Selmer structure.\comment{\footnote{In general, if one defines the dual 
Selmer structure for an arbitrary finite $R$-module $M$  which is also a $G_K$-module 
(here, $R$ is a complete, Noetherian and local ring of residue characteristic $p$), one needs to define it on the dual module $M^* = \Hom(M, R(1))$, where 
$R(1) = \Hom(R, \mu_{p^\infty})$, since the Tate local pairing is $\H^1(K_v, M) \times \H^1(K_v, M^*) \ra R(1)$. In our situation, $E[p^m]$ is self-dual 
because of the existence of the Weil pairing $e : E[p^m] \times E[p^m] \ra \mu_{p^m}$, so we will define it on $E[p^m]$. See~\cite[\S 2.3]{mazur-rubin:kolyvagin_systems} and~\cite{howard:heeg} for the general definition.}}}  
This will be the Selmer structure $\cF^*$ on $E[p^m]$ whose local conditions are the exact orthogonal 
complements of the local conditions of $\cF$ under the Tate local pairing 
$$
\langle\,, \rangle_v : \H^1(K_v, E[p^m]) \times \H^1(K_v, E[p^m]) \ra \Z/p^m\Z.  
$$  
Note that $\cF^*$ is a well-defined Selmer structure because the unramified local condition $\H^1_{\ur}(K_v, E[p^m])$ is self-dual for 
every place $v \nmid p$, for which the representation $E[p^m]$ is unramified (see~\cite[Thm.1.2.6]{milne:duality}).    

\vspace{0.1in}

\noindent \textit{5. Examples of Selmer structures.} The first Selmer structure which will be used in our argument is the standard 
\emph{Kummer Selmer structure} $\cF$ on $E[p^m]$. It is defined by the local condition $\H^1_{\cF}(K_v, E[p^m]) := \H^1_{\Kum}(K_v, E[p^m])$ for every 
$v$ (it is well-defined by Proposition~\ref{prop:ur-pts}). By using the compatibility of the Tate local duality with the Weil pairing, one shows 
that the structure $\cF$ is self-dual, i.e., the local conditions are self-orthogonal at each place with respect to the Tate local pairing 
$\langle \,,\rangle_v$.

We introduce a new set of Selmer structures which we refer to as \emph{stringent Kummer structures} and which are defined as follows:

\begin{definition}[stringent Kummer structures]
Let $q_1 q_2 \dots q_s \mid N$ be distinct prime each of which divides the conductor of $E$. 
Define a new Selmer structure $\cF_{\lceil q_1 \dots q_s \rceil }$ on $E[p^m]$ in the same way as $\cF$ except that for each places $v \mid q_i$ 
($i = 1, \dots, s$), the local condition at $v$ is replaced by the stringent Kummer condition $\H^1_{\Kum^0}(K_v, E[p^m])$. 
\end{definition}

\vspace{0.1in}

\indent Our main observation towards the refinement of Kolyvagin's results is the following  

\begin{proposition}\label{prop:refine}
For any stringent Kummer Selmer structure $\cF_{\lceil q_1 \dots q_s \rceil }$ and any $c \in \Lambda_m$ one has 
$$
\kappa_{c, m} \in \H^1_{\cF_{\lceil q_1\dots q_s \rceil }(c)}(K, E[p^m]). 
$$
The same holds for the class $\tilde{\kappa}_{c, m}$ for each $c$, for which it is defined. 
\end{proposition}

\begin{proof}
The cohomology classes $\kappa_{c, m}$ are known to satisfy the local conditions for the Selmer structure $\cF(c)$  
(see~\cite[Prop.6.2]{gross:kolyvagin},~\cite[Lem.4.3]{mccallum:kolyvagin} or~\cite[Lem.1.7.3]{howard:heeg}). 
Thus, we only need to check that $\loc_v(\kappa_{c, m}) \in \H^1_{\Kum^0}(K_v, E[p^m])$ for all places $v \mid N$.  

By Lemma~\ref{lemma:loc-ur-pts} the group $E^0(K_v^{\ur})$ is $p$-divisible, i.e., there is a short exact sequence 
$$
0 \ra E^0(K_v^{\ur})[p^m] \ra E^0(K_v^{\ur}) \xra{p^m} E^0(K_v^{\ur}) \ra 0. 
$$
By taking the long exact sequence on Galois cohomology and using the fact that N\'eron models are stable under \'etale base 
change, we obtain the following exact sequence 
$$
E^0(K_v) \ra \H^1(K_v^{\ur} / K_v, E^0(K_v^{\ur})[p^m]) \ra \H^1(K_v^{\ur} / K_v, E^0(K_v^{\ur}))[p^m] \ra 0. 
$$
But $\H^1(K_v^{\ur} / K_v, E^0(K_v^{\ur})) = 0$ according to Lang's theorem (see~\cite{lang:finitefields}), i.e., 
the map $E^0(K_v) \ra \H^1(K_v^{\ur} / K_v, E^0(K_v^{\ur})[p^m])$ is surjective. We thus consider the following 
commutative diagram 
$$
\xymatrix{
E^0(K_v) \ar[r] \ar[d] & \H^1(K_v^{\ur} / K_v, E^0(K_v^{\ur})[p^m]) \ar[d]^{\phi} & &  \\
E(K_v) \ar[r]& \H^1(K_v, E[p^m]) \ar[r] & \H^1(K_v, E)[p^m] \ar[r] & 0, \\     
} 
$$
where the map $\phi$ is the composition 
$$
\H^1(K_v^{\ur}/K_v, E^0(K_v^{\ur})[p^m]) \ra \H^1(K_v^{\ur}/K_v, E(K_v^{\ur})[p^m]) \xra{\inf} \H^1(K_v, E[p^m]).   
$$
We will be done if we show that the class $\loc_v(\kappa_{c, m}) \in \im (\phi)$ (the surjectivity then implies that it comes from a 
point of $E^0(K_v)$). 
To see that $\loc_v(\kappa_{c, m}) \in \im(\phi)$, we look at the explicit cocycle (\ref{eqn:cocycle}) and use Corollary~\ref{cor:red} 
according to which there exists a point 
$Q \in E(\Q)_{\tor}$, such that $Q_c = P_c - Q \in E^0(K_v^{\ur})$. Since $E(\Q)[p^\infty] = 0$ then the point $Q$ is $p$-divisible over $\Q$. This, together 
with the $p$-divisibility of $E^0(K_v^{\ur})$ (Lemma~\ref{lemma:loc-ur-pts}) implies  
$$
\kappa_{c, m}(\sigma) = -\frac{(\sigma - 1)P_c}{p^m} + \sigma \left ( \frac{P_c}{p^m} \right ) - \frac{P_c}{p^m} = -\frac{(\sigma - 1)Q_c}{p^m} + \sigma \left ( \frac{Q_c}{p^m} 
\right ) - \frac{Q_c}{p^m} \in E^0(K_v^{\ur})[p^m],     
$$
i.e., the class $\kappa_{c, m}$ is in the image (under $\phi$) of the cohomology class of $\H^1(K^{\ur}_v/K_v, E^0(K_v^{\ur})[p^m])$ represented by the 
cocycle $\ds \sigma \mapsto -\frac{(\sigma - 1)Q_c}{p^m} + \sigma \left ( \frac{Q_c}{p^m} \right ) - \frac{Q_c}{p^m}$. This proves the proposition. 
\end{proof}

\section{Poitou-Tate Global Duality}\label{sec:global-duality}
\comment{
We first state an useful consequence of the global duality theorem of Poitou and Tate which will be used extensively in the 
argument. Next, we derive an important consequence which allows us to compare the Selmer modules for two Selmer 
structures of the form $\cF(c)$ and $\cF(c\ell)$ where $c \in \Lambda$ and 
$\ell \in \Lambda^1$ (note that in the language of Mazur and Rubin the conductors $c$ and $c\ell$ are refered to as adjacent vertices in the Selmer 
graph which they consider).  
}

\subsection{Comparing Selmer modules and their duals}

The standard references for the Poitou-Tate global duality theorem are~\cite[Thm.1.7.3]{rubin:book}, \cite[Thm.I.4.10]{milne:duality} 
and~\cite[Thm.3.1]{tate:duality}. Here, we state a theorem which is an immediate consequence of global duality (see 
also~\cite[Thm.1.7.3]{rubin:book}). Before we stay the theorem, we introduce the following notation: if $\cF$ and $\cG$ are two Selmer structures,  
we say that $\cF \preceq \cG$ if $\cH^1_{\cF}(K_v, E[p^m]) \subseteq \cH^1_{\cG}(K_v, E[p^m])$ for every place $v$ of $K$. 

\begin{theorem}\label{thm:poitou-tate}
Let $\cF \preceq \cG$ be two Selmer structures on $E[p^m]$. 
Then there are exact sequences 
$$
0 \ra \H^1_{\cF}(K, E[p^m])^{\pm} \hra \H^1_{\cG}(K, E[p^m])^{\pm} \xra{(\loc^{\cG}_{\cF})^{\pm}} \bigoplus_{v} 
\left ( \frac{\H^1_{\cG}(K_v, E[p^m])}{\H^1_{\cF}(K_v, E[p^m])} \right )^{\pm}
$$
and 
$$
0 \ra \H^1_{\cG^*}(K, E[p^m])^{\pm} \hra \H^1_{\cF^*}(K, E[p^m])^{\pm} \xra{(\loc^{\cF^*}_{\cG^*})^{\pm}} \bigoplus_{v} 
\left ( \frac{\H^1_{\cF^*}(K_v, E[p^m])}{\H^1_{\cG^*}(K_v, E[p^m])} \right )^{\pm},  
$$ 
where $(\loc^{\cG}_{\cF})^{\pm}$ and 
$(\loc^{\cF^*}_{\cG^*})^{\pm}$ denote the natural restriction maps on the $\pm$-eigenspaces for complex conjugation
and the sum is over all places $v$, for which $\H^1_{\cF}(K_v, E[p^m]) \subsetneq \H^1_{\cG}(K_v, E[p^m])$. Moreover, the images of 
$(\loc^{\cG}_{\cF})^{\pm}$ and $(\loc^{\cF^*}_{\cG^*})^{\pm}$ are exact orthogonal complements with respect to the local pairings 
$\ds \sum_v \langle \,, \rangle^{\pm}_v$ obtained from the Tate pairings on the $\pm$-parts of the local cohomology groups.     
\end{theorem}

\subsection{Lozenge diagrams} 

The following result is a refinement of \cite[Lem.4.1.6]{mazur-rubin:kolyvagin_systems} and \cite[Lem.5.1.8]{howard:heeg}, and is very useful 
whenever one needs to compare the structures of two global Selmer modules for the Selmer structure $c$ and $c\ell$. 

\begin{lemma}\label{lem:lozenge}
Let $\cF$ be a Selmer structure on $E[p^m]$ (not necessarily self-dual). Consider the following diagrams
$$
\xymatrix{
& \H^1_{\cF^\ell(c)}(K, E[p^m])^{\pm} &  \\ 
\H^1_{\cF(c)}(K, E[p^m])^{\pm} \ar@{^{(}->}[ru]^{a^{\pm}} & & \ar@{_{(}->}[ul]_{b^{\pm}} \H^1_{\cF(c\ell)}(K, E[p^m])^{\pm} \\ 
& \H^1_{(\cF)_\ell(c)}(K, E[p^m])^{\pm} \ar@{_{(}->}[ul]_{c^{\pm}} \ar@{^{(}->}[ur]^{d^{\pm}} &       
} 
$$
and 
$$
\xymatrix{
& \H^1_{(\cF)_\ell(c)^*}(K, E[p^m])^{\pm} & \\ 
\H^1_{\cF(c)^*}(K, E[p^m])^{\pm} \ar@{^{(}->}[ru]^{(c^*)^{\pm}} & & \ar@{_{(}->}[ul]_{(d^*)^{\pm}} 
\H^1_{\cF(c\ell)^*}(K, E[p^m])^{\pm} \\
& \H^1_{\cF^\ell(c)^*}(K, E[p^m])^{\pm} \ar@{_{(}->}[ul]_{(a^*)^{\pm}} \ar@{^{(}->}[ur]^{(b^*)^{\pm}}, & 
}
$$
where each inclusion is labelled with the lengths of the corresponding cyclic cokernels. Then the lengths satisfy 

(i) $0 \leq a^{\pm}, b^{\pm}, c^{\pm}, d^{\pm}, (a^*)^{\pm}, (b^*)^{\pm}, (c^*)^{\pm}, (d^*)^{\pm} \leq m$.  

(ii) $a^{\pm} + c^{\pm} = b^{\pm} + d^{\pm}$ and $(a^*)^{\pm} + (c^*)^{\pm} = (b^*)^{\pm} + (d^*)^{\pm}$. 

(iii) $a^{\pm} + (a^*)^{\pm} = b^{\pm} + (b^*)^{\pm} = c^{\pm} + (c^*)^{\pm} = d^{\pm} + (d^*)^{\pm} = m$. 

(iv) $a^{\pm} \geq d^{\pm}$, $b^{\pm} \geq c^{\pm}$, $(c^*)^{\pm} \geq (b^*)^{\pm}$ and 
$(d^*)^{\pm} \geq (a^*)^{\pm}$. 
\end{lemma}

\begin{proof}
Statement (i) follows from the definition of $\cF_\ell(c)$, $\cF^{\ell}(c)$, $\cF^{\ell}(c)^*$ and $\cF_{\ell}(c)^*$ and the fact that the $\pm$-parts 
$\H^1_{\ur}(K_\lambda, E[p^m])^{\pm}$ and $\H^1_{\tr}(K_\lambda, E[p^m])^{\pm}$ of the unramified and the transverse local conditions are free 
of rank one over $\Z / p^m \Z$. Statement (ii) follows immediately from the diagram. Statement (iii) is an immediate consequence of the split global 
duality Theorem~\ref{thm:poitou-tate}. Finally, (iv) follows from the following equalities (which are consequences of the self-duality and the 
non-intersection of the transverse and the unramified local conditions) 
$$
\H^1_{\cF(c)}(K, E[p^m])^{\pm} \cap \H^1_{\cF(c\ell)}(K, E[p^m])^{\pm} = \H^1_{\cF_\ell(c)}(K, E[p^m])^{\pm}
$$
and 
$$
\H^1_{\cF(c)^*}(K, E[p^m])^{\pm} \cap \H^1_{\cF(c\ell)^*}(K, E[p^m])^{\pm} = \H^1_{\cF^\ell(c)^*}(K, E[p^m])^{\pm}. 
$$
\end{proof}

\section{The Main Theorem}\label{sec:main-thm}

In this section of the paper we prove Theorem~\ref{thm:main}. 

\subsection{An application of \v{C}ebotarev Density Theorem}

The following lemma is an application of \v{C}ebotarev density theorems which will be used in the proof of 
our theorem.  

\begin{lemma}\label{lem:cebotarev}
Assume Hypothesis (*) from the introduction and let  
$$
\kappa^{+} \in \H^1(K, E[p^m])^{+}, \  \kappa^{-} \in \H^1(K, E[p^m])^{-}
$$
be two cohomology classes with $M^{+} = \ord (\kappa^{+})$ and 
$M^{-} = \ord (\kappa^{-})$. Then there exists a prime $\ell \in \Lambda^1_{m}$, 
such that  $\ord ( \loc_\lambda(\kappa^{+}) ) = M^{+}$ and $\ord (\loc_\lambda(\kappa^{-})) = 
M^{-} \ne 0$.    
\end{lemma}

\begin{proof}
This follows immediately from \cite[Cor.~3.2]{mccallum:kolyvagin} since the cohomology classes $\kappa^+$ and $\kappa^-$ 
are linearly independent. 
\end{proof}

\begin{remark}
Notice that Hypothesis (*) can be weakened as the proof of~\cite[Cor.~3.2]{mccallum:kolyvagin} does not need the surjectivity of 
the Galois representation, but simply the weaker assumption that $\End_{\F_p}(E[p])$ is spanned (as an $\F_p$-vector space) by the 
elements $\sigma \in \Gal(\Q(E[p])/\Q)$. This fact is equivalent to the absolute irreducibility of the Galois representation 
$\rho_{E, p}$. 
\end{remark}

\subsection{Core vertices and minimal core vertices}
Let $m$ be an integer. For clarity, we denote each Selmer module $\H^1_{\cG(c)}(K, E[p^m])$ simply by 
$\cH_{\cG(c)}$ for various Selmer structures $\cG$ (i.e., we omit the Galois group and the Galois representation since they will stay fixed).

By a \emph{core vertex} for $m$ and the Kummer Selmer structure $\cF$ we mean any conductor $c \in \Lambda_m$, such that either 
$$
\Inv \cH_{\cF(c)}^{+} = (m) \  \textrm{and} \ \Inv \cH_{\cF(c)}^{-} = 0, 
$$ 
or 
$$ 
\Inv \cH_{\cF(c)}^{+} = (0) \  \textrm{and} \ \Inv \cH_{\cF(c)}^{-} = (m).    
$$
Moreover, a \emph{minimal core vertex} is a core vertex $c$ for $m$ for which $m(c) = m_\infty$ and $M(c) \geq m + m_\infty$. To make the argument 
easier to follow, we prove the theorem in the case when there exists a minimal core vertex for sufficiently large $m$. 

\begin{theorem}\label{thm:min-core}
Assume that there exists a minimal core vertex $c \in \Lambda_m$ for some $m > \max(m_{\max},  m_\infty)$. Then $m_\infty \geq m_{\max}$. 
\end{theorem}

\begin{proof}
Since $m > m(c) = m_\infty$, the class $\kappa_{c, m} \in \cH_{\cF(c)}^{\eps(c)}$ is non-trivial, i.e., we conclude that 
$\cH_{\cF(c)}^{\eps(c)} \isom \Z / p^m \Z$ and $\cH_{\cF(c)}^{-\eps(c)} = 0$.  
Consider the stringent Kummer Selmer structure $\cF_0 := \cF_{\lceil q \rceil}$ which differs from $\cF$ only at the two places 
$v$ and $\overline{v}$ above $q$ for which 
$\ord_p(c_q) = m_{\max}$ (note that the Tamagawa numbers $c_{v}$ and $c_{\overline{v}}$ are both equal to $c_q$). Consider the 
two Selmer modules $\cH_{\cF_0(c)}$ and $\cH^1_{\cF_0(c)^*}$. 

Since $\cH_{\cF_0(c)}^{\eps(c)}$ is a submodule of $\cH_{\cF(c)}^{\eps(c)}$ then  $\cH_{\cF_0(c)}^{\eps(c)} \isom \Z / p^{m'} \Z$ 
for some integer $m'$ satisfying $0 \leq m' \leq m$. Similarly, $\cH_{\cF_0(c)}^{-\eps(c)} = 0$.  
    
One can now determine the invariants of the dual modules $\cH_{\cF_0(c)^*}^{\pm}$ in terms of $m, m'$ and $m_{\max}$ by applying 
Theorem~\ref{thm:poitou-tate} to the exact sequence\footnote{Note that the exact sequence implies $m - m' \leq m_{\max}$.}
$$
0 \ra \cH_{\cF_0(c)}^{\pm \eps(c)} \hra \cH_{\cF(c)}^{\pm \eps(c)} \ra  \frac{\H^1_{\Kum}(K_v, E[p^m])}{\H^1_{\Kum^0}(K_v, E[p^m])}  
$$ 
and the dualized sequence 
$$
0 \ra \cH_{\cF(c)}^{\pm \eps(c)} \hra \cH_{\cF_0(c)^*}^{\pm \eps(c)} \ra \frac{\H^1_{\Kum^0}(K_v, E[p^m])^\perp}{\H^1_{\Kum}(K_v, E[p^m])},    
$$
where $\H^1_{\Kum^0}(K_v, E[p^m])^\perp$ is the orthogonal complement of $\H^1_{\Kum^0}(K_v, E[p^m])$ under the Tate local pairing. 
Since $p^{m_{\max}} \mid c_v$ and $m > m_{\max}$, we obtain (using the definition of the stringent Kummer Selmer structure and the fact that 
$c_v = [E(K_v) : E^0(K_v)]$) that both of the local quotients in the above sequences are isomorphic to $\Z / p^{m_{\max}}\Z$. This allows us to 
apply global duality (Theorem~\ref{thm:poitou-tate}) to conclude that 
$$
\Inv \left ( \cH_{\cF_0(c)^*}^{\eps(c)} \right ) = (m, m_{\max} + m' - m) \ \textrm{and} \ \Inv \left ( \cH_{\cF_0(c)^*}^{-\eps(c)} \right ) = (m_{\max}). 
$$
Moreover, the first invariant of $\cH_{\cF_0(c)^*}^{\eps(c)}$ corresponds to the free $\Z / p^m \Z$-submodule which contains 
$\kappa_{c, m}$. Since $m(c) + m \leq M(c)$, the class $\tilde{\kappa}_{c, m}$ is defined and generates the free submodule corresponding to 
the first invariant.  

Next, we use Lemma~\ref{lem:cebotarev} to choose a prime $\ell \in \Lambda^1_m$, such that $\ord( \loc_\lambda ( \tilde{\kappa}_{c, m}) ) = m$ and 
$\ord( \loc_{\lambda} ( \kappa' ) ) = m_{\max}$, where $\kappa'$ is a generator for the cyclic $\Z / p^m \Z$-module 
$\cH_{\cF_0(c)^*}^{-\eps(c)} \isom \Z / p^{m_{\max}} \Z$.

The key idea to finish the proof is to determine the invariants of the Selmer module 
$\cH_{(\cF_0)^\ell(c\ell)}^{\pm \eps(c)}$. Indeed, by the choice of $\ell$, the cyclic cokernels of the maps 
$\cH_{(\cF_0)^\ell(c)^*}^{\pm \eps(c)} \ra \cH_{\cF_0(c)^*}^{\pm \eps(c)}$ have lengths $m$ and $m_{\max}$, respectively 
(\comment{we say that the cokernel has type $(m; m_{\max})$}here, respectively means that the cokernel of the $\eps(c)$-map has length $m$ and 
the cokernel of the $-\eps(c)$-map has length $m_{\max}$).  
This means (by Lemma~\ref{lem:lozenge}(iii)) that the cokernels of the corresponding dual maps 
$\cH_{\cF_0(c)}^{\pm \eps(c)} \ra \cH_{(\cF_0(c))^\ell}^{\pm \eps(c)}$ have lengths 0 and $m - m_{\max}$, respectively. Thus, 
we conclude that $\cH_{(\cF_0)^\ell(c\ell)}^{- \eps(c)}$ is cyclic of length $m-m_{\max}$.  

\comment{
We summarize the computations of the cokernels of the various maps in the following three lozenge diagrams: 
$$
\xymatrix{
& \cH_{(\cF_0)_\ell(c)^*}^{\pm \eps(c)} 															\\ 
\cH_{\cF_0(c)^*}^{\pm \eps(c)}  \ar@{^{(}->}[ru]^{(m - m'; m)} & & \ar@{_{(}->}[ul]_{} 				
\cH_{\cF_0(c\ell)^*}^{\pm \eps(c)}   																\\
& \cH_{(\cF_0)^\ell(c)^*}^{\pm \eps(c)} \ar@{_{(}->}[ul]^{(m; m_{\max})} \ar@{^{(}->}[ur]^{}        \\ 
}
$$

$$
\xymatrix{
& \cH_{\cF^\ell(c)}^{\pm \eps(c)} & \\ 
\cH_{\cF(c)}^{\pm \eps(c)}\ar@{^{(}->}[ru]^{(0;m)} & & \ar@{_{(}->}[ul]_{} 
\cH_{\cF(c\ell)}^{\pm \eps(c)} \\
& \cH_{\cF_\ell(c)}^{\pm \eps(c)} \ar@{_{(}->}[ul]^{(m;0)} \ar@{^{(}->}[ur]^{} & 
}
$$

$$
\xymatrix{
& \cH_{(\cF_0)^\ell(c)}^{\pm \eps(c)}  & \\ 
\cH_{\cF_0(c)}^{\pm \eps(c)}  \ar@{^{(}->}[ru]^{(0;m-m_{\max})} & & \ar@{_{(}->}[ul]_{} 
\cH_{\cF_0(c\ell)}^{\pm \eps(c)} \\
& \cH_{(\cF_0)_\ell(c)}^{\pm \eps(c)} \ar@{_{(}->}[ul]^{(m';0)} \ar@{^{(}->}[ur]^{}, & 
}
$$
}

To complete the argument, notice $\cH^{-\eps(c)}_{\cF_0(c\ell)}$ is a submodule of 
$\cH^{-\eps(c)}_{(\cF_0)^{\ell}(c\ell)} \isom \Z / p^{m - m_{\max}}\Z$. Since $\kappa_{c\ell, m} \in \cH^{-\eps(c)}_{\cF_0(c\ell)}$ then 
$$
m - m_{\max} \geq \ord (\loc_{\lambda}(\kappa_{c\ell, m})) = \ord (\loc_{\lambda}(\kappa_{c, m})) = \ord (\kappa_{c, m}) = m - m_{\infty}, 
$$ 
where the second equality follows from Proposition~\ref{prop:comparison} and the third equality follows because 
$\ord (\tilde{\kappa}_{c, m}) = m$. Thus, $m_{\infty} \geq m_{\max}$. 
\end{proof}

\subsection{Existence of core vertices and the general case}\label{subsec:gencase}

In this final section, we reduce the proof of Theorem~\ref{thm:main} to the case when there exists a minimal core vertex (Theorem~\ref{thm:min-core}).  
Whenever $m$ is fixed, we accept the Selmer modules notation from the previous section. The most difficult part of the proof is the following technical 

\begin{proposition}\label{prop:core-vertex}
Let $c \in \Lambda$ satisfy $m(c) + m \leq M(c)$. 
There exists a core vertex $c' \in \Lambda_{m+m(c)}$, such that $m(c') \leq m(c)$.   
\end{proposition}

\begin{proof}
Since $m(c) + m \leq M(c)$, the class $\tilde{\kappa}_{c, m} \in \cH_{\cF(c)}$ generates a free $\Z / p^m \Z$-submodule of 
$\cH_{\cF(c)}$ which contains $\kappa_{c, m}$. This means that 
$$
\Inv \left (\cH_{\cF(c)}^{\eps(c)} \right ) = (m, x_1, x_2, \dots) 
$$  
for some $x_1 \geq x_2 \geq \dots \geq 0$. Let  
$$
\Inv \left (\cH_{\cF(c)}^{-\eps(c)} \right ) = (y_1, y_2, \dots), 
$$ 
where $y_1 \geq y_2 \geq \dots \geq 0$. Choose a prime $\ell_1 \in \Lambda^1_{M(c)}$, such that $\ord (\loc_{\lambda_1}(\kappa')) = y_1$ and $\ord(\loc_{\lambda_1}(\tilde{\kappa}_{c, m})) = m$ (here, $\lambda_1$ is the unique prime of $K$ above $\ell_1$ and $\kappa'$ is a generator 
for a direct summand of $\cH_{\cF(c)}^{-\eps(c)}$ corresponding to the invariant $y_1$). 
Such a prime can be selected according to Lemma~\ref{lem:cebotarev}. 

Next, we compute the invariants of $\cH_{\cF(c\ell_1)}^{\pm\eps(c\ell_1)}$. More precisely, we claim that $y_1 = y_2$ and  
$$
\Inv \left (\cH_{\cF(c\ell_1)}^{-\eps(c\ell_1)} \right ) = (x_3, x_4, \dots ) 
$$
and 
$$
\Inv \left (\cH_{\cF(c\ell_1)}^{\eps(c\ell_1)} \right ) = (m, y_1, y_2, \dots ).  
$$ 

To show this, we look at the lozenge diagrams for the self-dual Kummer Selmer structure $\cF$ and fill up as much as we can the lengths of 
the corresponding cokernels (the notation $(a; b)$ means that the cyclic cokernel for the $\eps(c)$-part has length $a$ and the 
cyclic cokernel for the $-\eps(c)$-part has length $b$ - we refer to it as the type of the cokernels)
$$
\xymatrix{
& \cH_{\cF^{\ell_1}(c)}^{\pm \eps(c)} & \\ 
\cH_{\cF(c)}^{\pm \eps(c)} \ar@{^{(}->}[ru]^{(0;m-y_1)} & & \ar@{_{(}->}[ul]_{(m; x)} 
\cH_{\cF(c\ell_1)}^{\pm \eps(c)} \\
& \cH_{\cF_{\ell_1}(c)}^{\pm \eps(c)} \ar@{_{(}->}[ul]^{(m;y_1)} \ar@{^{(}->}[ur]_{(0; m-x)} & 
}
$$ 
Indeed, the choice of $\ell_1$ implies that the cokernels of the maps $\cH_{\cF_{\ell_1}(c)}^{\pm \eps(c)} \hra \cH_{\cF(c)}^{\pm \eps(c)}$ 
have type $(m; y_1)$. By Lemma~\ref{lem:lozenge}(iii), the corresponding dual maps 
$\cH_{\cF(c)}^{\pm \eps(c)} \hra \cH_{\cF^{\ell_1}(c)}^{\pm \eps(c)}$ have cokernels of type $(0; m-y_1)$. Next, by Lemma~\ref{lem:lozenge}(iv), 
the cokernel of the map $\cH_{\cF(c\ell_1)}^{\eps(c)} \hra \cH_{\cF^{\ell_1}(c)}^{\eps(c)}$ has length $m$ and the cokernel of 
the map $\cH_{\cF_{\ell_1}(c)}^{\eps(c)} \hra \cH_{\cF(c)}^{\eps(c)}$ is trivial. Let $x$ be the length of the cokernel of the map 
$\cH_{\cF(c\ell_1)}^{-\eps(c)} \hra \cH_{\cF^{\ell_1}(c)}^{-\eps(c)}$. Again, 
by Lemma~\ref{lem:lozenge}(iii), the cokernel of the map $\cH_{\cF_{\ell_1}(c)}^{-\eps(c)} \hra \cH_{\cF(c)}^{-\eps(c)}$ has length $m-x$. 
This justifies the labelling of the diagram.   

The next crucial observation is that $x = y_1$. To see this, we observe that the comparison isomorphism $\phi_{\lambda_1}$ and 
Proposition~\ref{prop:comparison} imply
$$
m - m(c\ell_1) \geq \ord (\loc_{\lambda_1}(\kappa_{c\ell_1, m})) = \ord ( \loc_{\lambda_1}(\kappa_{c\ell_1, m})) = \ord ( \kappa_{c\ell_1, m}) = m - m(c),   
$$ 
i.e., $m(c\ell_1) \leq m(c)$. This means that $m(c\ell_1) + m \leq M(c) = M(c\ell_1)$, so the class $\kappa_{c\ell_1, m}$ lives in a free 
submodule of $\cH_{\cF(c\ell_1)}^{\eps(c\ell_1)}$ of rank one. Thus, the module $\cH_{\cF(c\ell_1)}$ contains an invariant $m$. Hence, 
by looking at the inclusion $\cH_{\cF_{\ell_1}(c)} \hra \cH_{\cF(c\ell_1)}$ and using that $x \geq y_1$, we obtain $x = y_1 = y_2$. Moreover, 
we determine the invariants 
$$
\Inv \cH_{\cF(c\ell_1)}^{\eps(c\ell_1)} = (m, y_3, \dots) \textrm{ and } \Inv \cH_{\cF(c\ell_1)}^{-\eps(c\ell_1)} = (x_1, x_2, \dots). 
$$  

We repeat the above process  and use Lemma~\ref{lem:cebotarev} to choose a prime $\ell_2 \in \Lambda_{M(c)}$, such that 
$\ord (\loc_{\lambda_2}(\kappa_2)) = x_1$ and $\ord(\loc_{\lambda_2}(\tilde{\kappa}_{c, m})) = m$, where $\lambda_2$ is the unique prime of 
$K$ lying above $\ell_2$ and $\kappa_2$ is a generator of 
a direct summand corresponding to the invariant $x_2$. By exactly the same argument, $m(c\ell_1 \ell_2) \leq m(c\ell_1) \leq m(c)$ and 
$$
\Inv \cH_{\cF(c\ell_1 \ell_2)}^{\eps(c\ell_1 \ell_2)} = (m, x_3, x_4, \dots) \textrm{ and } \Inv \cH_{\cF(c\ell_1\ell_2)}^{-\eps(c\ell_1 \ell_2)} = 
(y_3, y_4, \dots). 
$$
  
We continue the process of adding primes from $\Lambda_{M(c)}^1$ to the conductor until we reach at a conductor $c' = c \ell_1 \dots \ell_s$, such that 
$\cH_{\cF(c')}^{\eps(c')} \isom \Z / p^m \Z$ and $\cH_{\cF(c')}^{\eps(c')} = 0$ and for which $m(c') \leq m(c)$.  
\end{proof}

Finally, we prove the main Theorem~\ref{thm:main}. 

\begin{proof}[Proof of Theorem~\ref{thm:main}]
According to~\cite[Thm.1]{kolyvagin:structureofsha}, one can also describe $m_\infty$ as 
$$
m_\infty = \lim_{r \ra \infty} \inf_{c \in \Lambda_{m'}^r} m(c)  
$$
for any sufficiently large integer $m'$ (in particular, for $m' > 2m_0$). Fix an integer $m > \max\{m_0, m_{\max}\}$ and let $m' = m + m_0$. 
Choose $c \in \Lambda_{m'}^r$ for which $m(c) = m_\infty$ (such a $c$ exists according to the above redefinition of $m_\infty$). 

Since $m(c) + m = m_\infty + m < m' \leq M(c)$, by Theorem~\ref{prop:core-vertex}, there exists 
a core vertex $c' \in \Lambda_{m+m(c)}$, such that $m(c') \leq m(c)  = m_\infty$ (i.e., $m(c') = m_\infty$). This, together with $m(c') + m \leq M(c')$ 
implies that $c'$ is a minimal core vertex. Thus, Theorem~\ref{thm:min-core} implies the desired inequality.   
\end{proof}

\bibliography{biblio}
\end{document}

%% file: macros.tex
\usepackage[active]{srcltx}\usepackage[active]{srcltx}\usepackage[active]{srcltx}\usepackage[active]{srcltx}
\usepackage[centertags]{amsmath}
\usepackage{amsfonts}
\usepackage{amssymb}
\usepackage{amsthm}
\usepackage{amscd}

\DeclareMathOperator{\Inv}{\mathbf{Inv}}

\DeclareMathOperator{\loc}{loc}

\DeclareMathOperator{\Kum}{Kum}
\DeclareMathOperator{\f}{f}
\DeclareMathOperator{\s}{s}

\DeclareMathOperator{\chr}{char}

\DeclareMathOperator{\im}{im}

\DeclareMathOperator{\End}{End}

\DeclareMathOperator{\Pic}{Pic}

\DeclareMathOperator{\Sel}{Sel}

\DeclareFontEncoding{OT2}{}{} 
  \newcommand{\textcyr}[1]{%
    {\fontencoding{OT2}\fontfamily{wncyr}\fontseries{m}\fontshape{n}%
     \selectfont #1}}
\newcommand{\Sha}{{\mbox{\textcyr{Sh}}}}

\newcommand{\cH}{\mathcal{H}}

\newcommand{\cF}{\mathcal{F}}
\newcommand{\cG}{\mathcal{G}}

\newcommand{\cO}{\mathcal{O}}

\newcommand{\ds}{\displaystyle}

\DeclareMathOperator{\tor}{tor}

\DeclareMathOperator{\Jac}{Jac}
\DeclareMathOperator{\res}{res}
\DeclareMathOperator{\Ker}{Ker}

\newcommand{\Kbar}{\overline{K}}
\newcommand{\Lbar}{\overline{L}}

\newcommand{\ra}{\rightarrow}
\newcommand{\xra}[1]{\xrightarrow{#1}}
\newcommand{\hra}{\hookrightarrow}

\newcommand{\Gm}{\mathbb{G}_m}

\newcommand{\eps}{\varepsilon}
\newcommand{\vphi}{\varphi}

\newcommand{\comment}[1]{}
\newcommand{\Q}{\mathbb{Q}}

\newcommand{\cN}{\mathcal{N}}

\newcommand{\C}{\mathbb{C}}

\newcommand{\Qbar}{\overline{\Q}}
\newcommand{\kbar}{\overline{k}}

\newcommand{\Z}{\mathbb{Z}}
\newcommand{\F}{\mathbb{F}}

\newcommand{\isom}{\cong}
\DeclareMathOperator{\ab}{ab}

\DeclareMathOperator{\Fr}{Fr}

\DeclareMathOperator{\ord}{ord}
\DeclareMathOperator{\GL}{GL}

\DeclareMathOperator{\Gal}{Gal}

\newcommand{\m}{\mathfrak{m}}

\DeclareMathOperator{\tr}{tr}

\DeclareMathOperator{\ur}{ur}
\DeclareMathOperator{\Tr}{Tr}
\DeclareMathOperator{\Hom}{Hom}

\DeclareMathOperator{\HH}{H}
\renewcommand{\H}{\HH}

\theoremstyle{plain}
\newtheorem{theorem}{Theorem}[section]
\newtheorem{proposition}[theorem]{Proposition}
\newtheorem{corollary}[theorem]{Corollary}

\newtheorem{lemma}[theorem]{Lemma}

\newtheorem{conjecture}[theorem]{Conjecture}

\theoremstyle{definition}
\newtheorem{definition}[theorem]{Definition}

\newtheorem{alg}[theorem]{Algorithm}

\newtheorem{remark}[theorem]{Remark}

\numberwithin{equation}{section}
\numberwithin{figure}{section}
\numberwithin{table}{section}

\newcounter{listnum}

{\end{alg}}
{\begin{enumerate}\setlength{\itemsep}{0.1ex}}{\end{enumerate}}

